%% file: Entekhabi_Lancaster_Paper_2.tex
\documentclass{article}
\usepackage{graphicx}
\usepackage{amsmath}
\usepackage{amssymb}
\usepackage{amsthm}
\usepackage{epsfig}

\title{Radial Limits of Capillary Surfaces at Corners}
\author{Mozhgan (Nora) Entekhabi \& Kirk E. Lancaster      \\
                       Department of Mathematics, Statistics \& Physics \\
                            Wichita State University \\
                            Wichita, Kansas, 67260-0033}

\def\Real{{\rm I\hspace{-0.2em}R}}

\newcommand\myeq{\mathrel{\overset{\makebox[0pt]{\mbox{\normalfont\tiny\sffamily def}}}{=}}}

\date{ }

\begin{document}
\maketitle

\vspace*{3mm}


\vspace*{3mm}

\begin{abstract}
\noindent Consider a solution $f\in C^{2}(\Omega)$  of a prescribed mean curvature equation 
\[
{\rm div}\left(\frac{\nabla f}{\sqrt{1+|\nabla f|^{2}}}\right)=2H(x,f) \ \ \ \ {\rm in} \ \ \Omega\subset \Real^{2},
\]
where $\Omega$  is a domain whose boundary has a corner at ${\cal O}=(0,0)\in\partial\Omega$  
and the angular measure of this corner is $2\alpha,$  for some $\alpha\in (0,\pi).$  
Suppose $\sup_{x\in\Omega} |f(x)|$  and $\sup_{x\in\Omega} |H(x,f(x))|$  are both finite.  
If $\alpha>\frac{\pi}{2},$   then the (nontangential) radial limits of $f$  at ${\cal O},$  
\[
Rf(\theta) \myeq \lim_{r\downarrow 0} f(r\cos(\theta),r\sin(\theta)), 
\]
were recently proven by the authors to exist, independent of the boundary behavior of $f$  on $\partial\Omega,$   
and to have a specific type of behavior.   

Suppose $\alpha\in \left(\frac{\pi}{4},\frac{\pi}{2}\right],$  the contact angle $\gamma(\cdot)$  
that the graph of $f$  makes with one side of $\partial\Omega$  has a limit (denoted $\gamma_{2}$) at ${\cal O}$  
and 
\[
\pi-2\alpha < \gamma_{2} <2\alpha.
\]
We prove that the (nontangential) radial limits of $f$  at ${\cal O}$  exist 
and the radial limits have a specific type of behavior, 
independent of the boundary behavior of $f$  on the other side of $\partial\Omega.$ 
We also discuss the case $\alpha\in \left(0,\frac{\pi}{2}\right].$  
\end{abstract}

\newtheorem{thm}{Theorem}
\newtheorem{prop}{Proposition}
\newtheorem{cor}{Corollary}
\newtheorem{lem}{Lemma}
 
\section{Introduction and Statement of Main Theorems}

Let $\Omega$  be a domain  in ${\Real}^{2}$  whose boundary has a corner at ${\cal O}\in\partial\Omega.$  
Suppose $H:\Omega\times\Real \to \Real$  and $H$  satisfies one of the conditions which guarantees that ``cusp solutions'' 
(e.g. \S 5 of \cite{LS1}, \cite{LS2}) do not exist; for example,  $H({\bf x},t)$  is strictly increasing in $t$  
for each ${\bf x}$  or is real-analytic (e.g. constant).  We will assume ${\cal O}=(0,0).$  
Let $\Omega^{*} = \Omega \cap B_{\delta^{*}}({\cal O})$, where $B_{\delta^{*}}({\cal O})$ is the ball in $\Real^{2}$ 
of radius $\delta^{*}$ about ${\cal O}$.   
Polar coordinates relative to ${\cal O}$ will be denoted by $r$ and $\theta$.
We assume that $\partial \Omega$  is piecewise smooth and there exists $\alpha \in \left(0,\pi\right)$  
such that $\partial \Omega\setminus \{{\cal O}\} \cap B_{\delta^{*}}({\cal O})$  consists of two (open) $C^{1}$ arcs  
${\partial}^{+}\Omega^{*}$  
and $\partial^{-}\Omega^{*}$, whose tangent lines approach the lines $L^{+}: \:  \theta = \alpha$  and 
$L^{-}: \: \theta = - \alpha$, respectively, as the point ${\cal O}$ is approached.

Suppose $\alpha>\frac{\pi}{2},$  $f\in C^{2}(\Omega)$  satisfies the prescribed mean curvature equation
\begin{equation}
\label{PMC}
Nf(x)  =  2H(x,f(x))  \mbox{  \ for \ } x\in \Omega,    
\end{equation}
where $Nf = \nabla \cdot Tf = {\rm div}\left(Tf\right)$  and $Tf = \frac{\nabla f}{\sqrt{1 + |\nabla f|^{2}}},$   and 
\begin{equation}
\label{Bounds}
\sup_{x\in\Omega} |f(x)| < \infty \ \ \ \  {\rm and} \ \ \ \  \sup_{x\in\Omega} |H(x,f(x))| < \infty.  
\end{equation}
In \cite{NoraKirk1}, the authors proved that the radial limits
\[
Rf(\theta) \myeq \lim_{r\downarrow 0} f(r\cos(\theta),r\sin(\theta)) 
\]
exist for all $\theta\in (-\alpha,\alpha),$   $Rf(\cdot)$  is a continuous function on $(-\alpha, \alpha)$ 
and these radial limits have similar behavior to that observed in Theorem 1 of \cite{LS1}.

\begin{figure}[ht]
\centering
\includegraphics{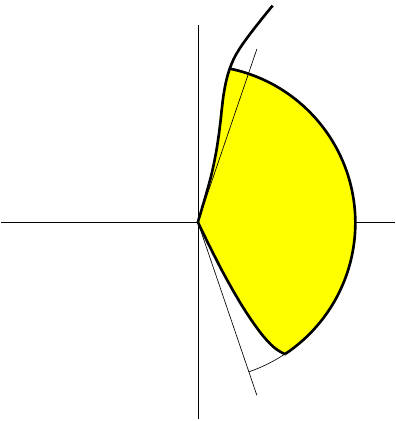}
\caption{The domain $\Omega^{*}$  \label{ZERO}}
\end{figure}

Suppose $\alpha\le \frac{\pi}{2}$  (see Figure \ref{ZERO}) and $f\in C^{2}(\Omega)\cap C^{1}\left(\Omega\cup \partial^{-}\Omega^{*}\right)$  
satisfies (\ref{PMC}) and (\ref{Bounds}). In \cite{NoraKirk1}, it is shown that if 
\begin{equation}
\label{Limit_of_f-}
\lim_{\partial^{-}\Omega^{*}\ni {\bf x}\to {\cal O} } f\left({\bf x}\right)=z_{2} \ \ \ \ {\rm exists},
\end{equation}
then the radial limits of $f$  at ${\cal O}$  exist and behave as expected.  In this paper, we consider the capillary problem as 
our model and suppose $f\in C^{2}(\Omega)\cap C^{1}\left(\Omega\cup \partial^{-}\Omega^{*}\right)$  
satisfies (\ref{PMC}), (\ref{Bounds}) and the boundary condition 
\begin{equation}
\label{PMCBC}
Tf(x)\cdot \nu(x)= \cos(\gamma(x))  \mbox{\ for \ } x\in \partial^{-}\Omega^{*},
\end{equation}
where  $\nu(x)$  is the exterior unit normal to $\Omega$  at $x\in  \partial\Omega$  and 
$\gamma:\partial\Omega\to [0,\pi]$  is the contact angle between the graph of $f$  and $\partial\Omega\times\Real,$   and 
\begin{equation}
\label{Limit-}
\lim_{\partial^{-}\Omega^{*}\ni {\bf x}\to {\cal O} } \gamma\left({\bf x}\right)=\gamma_{2}.
\end{equation}
We shall prove 

\begin{thm}
\label{CONCLUSION}
Let $f\in C^{2}(\Omega)\cap C^{1}\left(\Omega\cup \partial^{-}\Omega^{*}\right)$  
satisfy (\ref{PMC}) \& (\ref{PMCBC}) and suppose (\ref{Bounds}) and (\ref{Limit-}) hold, 
$\alpha \in \left(\frac{\pi}{4},\frac{\pi}{2}\right]$  and 
\begin{equation}
\label{Condition-1} 
\pi-2\alpha<\gamma_{2}<2\alpha.  
\end{equation}
Then (\ref{Limit_of_f-}) holds, $Rf(\theta)$  exists for all $\theta\in (-\alpha,\alpha)$  and 
$Rf(\cdot)$  is a continuous function on $[-\alpha, \alpha),$  
where $Rf(-\alpha)\myeq z_{2}.$   Further $Rf(\cdot)$  behaves in one of the following ways:

\noindent (i)  $Rf:[-\alpha,\alpha)\to\Real$  is a constant function (hence $f$  has a nontangential limit at ${\cal O}$). 

\noindent (ii) There exist $\alpha_{1}$ and $\alpha_{2}$ so that $-\alpha \leq \alpha_{1}
< \alpha_{2} \leq \alpha$ and $Rf$ is constant on $[-\alpha, \alpha_{1}]$ and
$[\alpha_{2}, \alpha)$ and strictly increasing or strictly decreasing on
$[\alpha_{1}, \alpha_{2})$.  
\end{thm}

\noindent If $\alpha \in \left(0,\frac{\pi}{4}\right],$  then (\ref{Condition-1}) cannot be satisfied. 
If $\alpha \in \left(\frac{\pi}{4},\frac{\pi}{2}\right]$  but $\gamma_{2}\ge 2\alpha$  or $\gamma_{2}\le \pi-2\alpha,$
then (\ref{Condition-1}) is not satisfied.  In both cases, Theorem \ref{CONCLUSION} is not applicable.
In these cases, we can prove the existence of $Rf(\cdot)$ if we add an assumption about the behavior 
of $\gamma$  on $\partial^{+}\Omega^{*}.$

\begin{thm}
\label{DEMONSTRATION}
Let $f\in C^{2}(\Omega)\cap C^{1}\left(\Omega\cup \partial^{-}\Omega^{*}\cup \partial^{+}\Omega^{*}\right)$  
satisfy (\ref{PMC})-(\ref{PMCBC}).  Suppose (\ref{Bounds}) and (\ref{Limit-}) hold, 
$\alpha \in \left(0,\frac{\pi}{2}\right],$  there exist $\lambda_{1},\lambda_{2}\in [0,\pi]$  
with $0<\lambda_{2}-\lambda_{1}<4\alpha$  such that $\lambda_{1}\le \gamma(x)\le \lambda_{2}$  
for $x\in \partial^{+}\Omega^{*}$  and 
\begin{equation}
\label{Condition-2} 
\pi-2\alpha-\lambda_{1}<\gamma_{2}<\pi+2\alpha-\lambda_{2}.  
\end{equation}
Then the conclusions of Theorem \ref{CONCLUSION} hold.
\end{thm}
\vspace{3mm}

\noindent {\bf Remark:} {\it Theorem \ref{DEMONSTRATION} only offers a new result when $\lambda_{1}=0$  or $\lambda_{2}=\pi;$  
Figure 8 of \cite{Shi} illustrates one example in which  $\lambda_{1}=0$  or $\lambda_{2}=\pi$  occurs. 
If $0<\lambda_{1}<\lambda_{2}<\pi,$  then Theorem \ref{DEMONSTRATION} is a 
consequence of Theorem 1 of \cite{LS1}; in this case, the argument here or in \cite{LS1} implies 
$Rf(\theta)$  exists for all $\theta\in [-\alpha,\alpha].$ }
\vspace{3mm}

\noindent {\bf Remark:} {\it In \cite{CF:94,Finn:96}, Concus and Finn proved that, in a neighborhood ${\cal U}$  of ${\cal O}$  
and assuming $\partial^{+}\Omega^{*}$  and $\partial^{-}\Omega^{*}$  are straight line segments, 
a solution of a constant mean curvature equation (i.e. $H$  is constant in (\ref{PMC})) with constant contact angles 
$\gamma_{1}$  on  ${\cal U}\cap \partial^{+}\Omega^{*}$  and $\gamma_{2}$  on  ${\cal U}\cap \partial^{-}\Omega^{*}$  
can exist only if $|\pi-\gamma_{1}-\gamma_{2}|\le 2\alpha.$    
Using this, when $\gamma_{1}=0,$   we would 
obtain a (local) upper bound for $f$  in Theorem \ref{CONCLUSION} when  $\pi-2\alpha< \gamma_{2}$  and, 
when $\gamma_{1}=\pi,$  a (local) lower bound for $f$ when $\gamma_{2}< 2\alpha;$  these two inequalities are 
equivalent to (\ref{Condition-1}). }
\vspace{3mm}

\noindent {\bf Remark:} {\it  As in \cite{LS1}, conclusion (\ref{Limit_of_f-}) of Theorems \ref{CONCLUSION} and \ref{DEMONSTRATION} 
is a consequence of a general argument; establishing (\ref{Limit_of_f-}) is not a key step in the proof.  }
\vspace{3mm}

\noindent {\bf Remark:} {\it One might contemplate replacing hypothesis (\ref{Limit-}) by something like 
$0<\sigma_{1}\le \gamma(x)\le \sigma_{2}<\pi$  for $x\in \partial^{+}\Omega^{*}$  (as in \cite{LS1}) and suitably 
modifying (\ref{Condition-1}) or (\ref{Condition-2}).  The comparison methods used here allow for this possibility.}

\section{Preliminary Remarks}

Let $f\in C^{2}(\Omega)$  satisfy  (\ref{PMC}) and suppose (\ref{Bounds})  holds.  
Throughout the remainder of the article, let us assume that $M_{1}\in (0,\infty),$    $M_{2}\in [0,\infty),$   
\begin{equation}
\label{NewBounds}
\sup_{x\in\Omega} |f(x)| \le M_{1} \ \ \ \  {\rm and} \ \ \ \  \sup_{x\in\Omega} |H(x,f(x))| \le M_{2}.  
\end{equation}

\subsection{A Specific Torus}

We will use portions of tori and comparison function arguments as, for instance, in Examples 2 \& 3 of \cite{LS1} 
and the Courant-Lebegsue lemma (\cite{Cour:50}, Lemma 3.1) to obtain upper and lower bounds on $f$  near ${\cal O}$  in 
specific subsets of $\Omega$  and prove Theorems \ref{CONCLUSION} and \ref{DEMONSTRATION}.    
Let us discuss the construction of a particular torus. 

\begin{figure}[ht]
\centering
\includegraphics{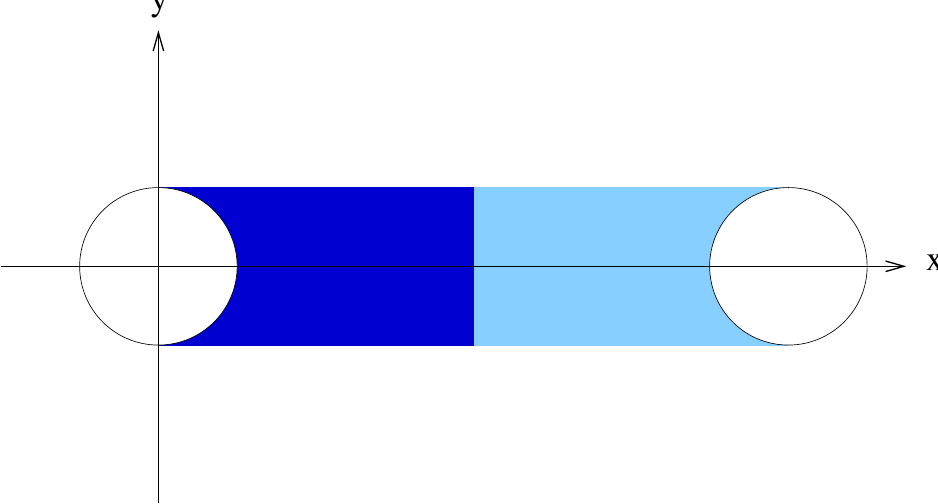}
\caption{The regions $\Delta$  (in dark blue) and $\Delta^{R}$  (in light blue) \label{ONEa}}
\end{figure}

\begin{figure}[ht]
\centering
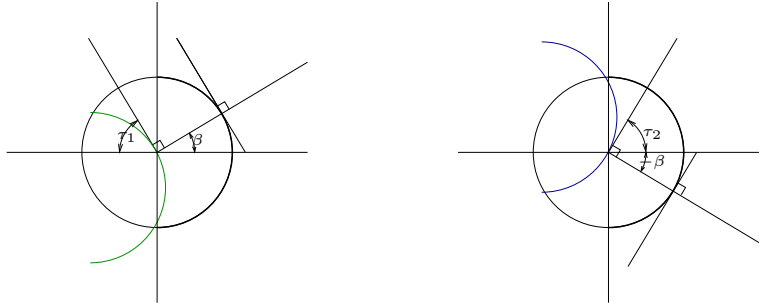
\caption{Left:  $\beta+\tau_{1}=\frac{\pi}{2}$ \hspace{1cm} Right: $-\beta+\tau_{2}=\frac{\pi}{2}$  \hspace{1cm} ($\beta\ge 0$)  
\label{ZEROX0}}
\end{figure}

\noindent Set 
\[
r_{0}=\left\{ \begin{array}{cl}
1
& \mbox{{\rm if}\ $M_{2}=0$  }  \\ 
  &  \\ 
\frac{1}{M_{2}}+1-\sqrt{\left(\frac{1}{M_{2}}\right)^{2}+1}
& \mbox{{\rm if}\ $M_{2}>0.$}  
\end{array}   \right.  \]
Let $\Delta=\{{\bf x}=(x_{1},x_{2})\in\Real^{2} \ : \ |{\bf x}|\ge r_{0}, \ 0\le x_{1}\le 2, \ |x_{2}|\le r_{0} \}$
and $\Delta^{R}=\{{\bf x}=(x_{1},x_{2})\in\Real^{2} \ : \ (4-x_{1},x_{2})\in\Delta\}.$
Let 
\[
{\cal T}= \{ (2+(2+r_{0}\cos(v))\cos(u),r_{0}\sin(v),(2+r_{0}\cos(v))\sin(u))\ : \  u,v\in [0,2\pi] \}
\]
be a torus with axis of symmetry $\{(2,y,0) : y\in\Real\},$  major radius $R_{0}=2$  and minor radius $r_{0};$  
recall that the mean curvature of ${\cal T}$  (with respect to the exterior normal) 
at $(2+(2+r_{0}\cos(v))\cos(u),r_{0}\sin(v),(2+r_{0}\cos(v))\sin(u))$  is 
\[
\frac{1}{2}H_{T}=-\frac{2+2r_{0}\cos(v)}{2r_{0}(2+r_{0}\cos(v))}.
\]
A calculation shows that  
\begin{equation}
\label{Limits}
-\left(\frac{1}{r_{0}}+\frac{1}{2+r_{0}}\right) \le H_{T} \le -\left(\frac{1}{r_{0}}-\frac{1}{2-r_{0}}\right) = -M_{2}.
\end{equation}
Set 
\[
{\cal T}^{+}= \{({\bf x},z)\in {\cal T} \ : \ {\bf x}\in\Delta, \ z\ge 0\} \ \ \ {\rm and} \ \ \ 
{\cal T}^{-}= \{({\bf x},z)\in {\cal T}  \ : \ {\bf x}\in\Delta, \ z\le 0\}.
\]
Let $h^{+},h^{-}:\Delta \to \Real$  be functions whose graphs satisfy 
\[
\{({\bf x},h^{+}({\bf x})) : {\bf x}\in\Delta \} = {\cal T}^{+} \ \ \  {\rm and} \ \ \ 
\{({\bf x},h^{-}({\bf x})) : {\bf x}\in\Delta \} = {\cal T}^{-}. 
\]
Then, from (\ref{Limits}), we have  
\begin{equation}
\label{MeanBounds}
{\rm div}\left(\frac{h^{+}}{\sqrt{1+|\nabla h^{+}|^{2}}}\right) \ge M_{2}
\ \ \ \  {\rm and} \ \ \ \    
{\rm div}\left(\frac{h^{-}}{\sqrt{1+|\nabla h^{-}|^{2}}}\right) \le -M_{2}.
\end{equation}
For each $\beta\in \left[-\frac{\pi}{2},\frac{\pi}{2}\right]$  let 
$\Delta_{\beta} = {\cal R}_{\alpha} \circ T_{\beta}\left(\Delta\right),$  
where ${\cal R}_{\alpha}:\Real^{2}\to \Real^{2},$  given by 
\[
(x_{1},x_{2})\mapsto (\cos(\alpha)x_{1}+\sin(\alpha)x_{2},-\sin(\alpha)x_{1}+\cos(\alpha)x_{2}),
\]
is the rotation about $(0,0)$  through the angle $-\alpha$  and $T_{\beta}:\Real^{2}\to \Real^{2},$  given by 
\[
(x_{1},x_{2})\mapsto (x_{1}-r_{0}\cos(\beta),x_{2}-r_{0}\sin(\beta)),
\]
is the translation taking $(r_{0}\cos(\beta),r_{0}\sin(\beta))\in \partial\Delta$  to $(0,0).$  
We will let $\tau_{1}$  denote the angle that upward tangent ray to $T_{\beta}(C)$  makes with the negative 
$x_{1}-$axis and let $\tau_{2}$  denote the angle that upward tangent ray to $T_{-\beta}(C)$  makes with the 
positive $x_{1}-$axis, where 
$C=\{{\bf x}=(x_{1},x_{2})\in\Real^{2} \ : \ |{\bf x}|= r_{0}, \ x_{1}\ge 0\}.$
(Figure \ref{ZEROX0} illustrates this when $\beta>0.$)  
Let $h^{\pm}_{\beta}:\Delta_{\beta}\to\Real$  be defined by $h^{\pm}_{\beta}=h^{\pm}\circ T_{\beta}^{-1}\circ {\cal R}_{\alpha}^{-1}.$  

Let $q$  denote the modulus of continuity of $h^{-}$  
(i.e. $|h_{\beta}^{-}({\bf x}_{1}) - h_{\beta}^{-}({\bf x}_{2})| \leq q(| {\bf x}_{1} - {\bf  x}_{2}|).$  
Notice that $q$  is also the modulus of continuity of $h^{+},$  as well as for $h_{\beta}^{-}$  and $h_{\beta}^{+}$
for each $\beta\in \left[-\frac{\pi}{2},\frac{\pi}{2}\right].$

\begin{figure}[ht]
\centering
\includegraphics{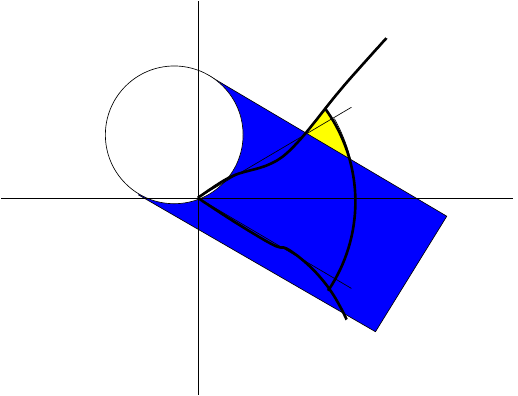}
\caption{The domain (in blue) of a toroidal function $h^{\pm}_{\beta},$   $\alpha<\frac{\pi}{4}.$  \label{ONEb}}
\end{figure}

\subsection{Parametric Framework}

Since $f\in C^{0}(\Omega),$   we may assume that $f$  is uniformly continuous on 
$\{{\bf x}\in \Omega^{*} \ :\ |{\bf x}|>\delta\}$  for each $\delta\in (0,\delta^{*});$  
if this is not true, we may replace $\Omega$  with $U,$  $U\subset \Omega,$  such that 
$\partial\Omega \cap \partial U = \{ {\cal O} \}$  and $\partial U \cap B_{\delta^{*}}({\cal O})$  
consists of two arcs  ${\partial}^{+} U$  and $\partial^{-} U$, whose tangent lines approach the lines 
$L^{+}: \:  \theta = \alpha$  and $L^{-}: \: \theta = - \alpha$, respectively, as the point ${\cal O}$ is approached.
Set
\[ 
S^{*}_{0} = \{ ({\bf x},f({\bf x})) : {\bf x} \in \Omega^{*} \}
\]
and
\[ 
\Gamma^{*}_{0} = \{ ({\bf x},f({\bf x})): {\bf x} \in \partial \Omega^{*}    \setminus \{ {\cal O} \} \}; 
\]
the points where $\partial B_{\delta^{*}}({\cal O})$ intersect $\partial \Omega$  are labeled $A\in {\partial}^{-}\Omega^{*}$ and 
$B \in {\partial}^{+}\Omega^{*}.$
From the calculation on page 170 of \cite{LS1},  we see that the area of $S^{*}_{0}$  is finite; let $M_{0}$  denote this area. 
For $\delta\in (0,1),$  set 
\[
p(\delta) = \sqrt{\frac{8\pi M_{0}}{\ln\left(\frac{1}{\delta}\right)}}.
\]
Let $E= \{ (u,v) : u^{2}+v^{2}<1 \}.$ 
As in \cite{EL:86,LS1}, there is a parametric description of the surface $S^{*}_{0},$  
\begin{equation}
\label{PARAMETRIC}
Y(u,v) = (a(u,v),b(u,v),c(u,v)) \in C^{2}(E:{\Real}^{3}), 
\end{equation}
which has the following properties:
 
\noindent $\left(a_{1}\right)$  $Y$ is a diffeomorphism of $E$ onto $S^{*}_{0}$.
 
\noindent $\left(a_{2}\right)$
Set $G(u,v)=(a(u,v),b(u,v)),$  $(u,v)\in E.$   Then $G \in C^{0}(\overline{E} : {\Real}^{2}).$  

\noindent $\left(a_{3}\right)$
Let $\sigma=G^{-1}\left(\partial \Omega^{*}\setminus \{ {\cal O} \}\right);$  
then $\sigma$ is a connected arc of $\partial E$  and $Y$ maps $\sigma$ strictly monotonically onto 
$\Gamma^{*}_{0}.$  
We may assume the endpoints of $\sigma$  are ${\bf o}_{1}$  and ${\bf o}_{2}$  and there exist points 
${\bf a}, {\bf b}\in\sigma$  such that $G({\bf a})=A,$  $G({\bf b})=B,$
$G$  maps the (open) arc ${\bf o}_{1}{\bf b}$  onto $\partial^{+}\Omega,$  and $G$ 
maps the (open) arc ${\bf o}_{2}{\bf a}$  onto $\partial^{-}\Omega.$ 
(Note that ${\bf o}_{1}$  and ${\bf o}_{2}$  are not assumed to be distinct at this point; 
one of Figure 4a or 4b of \cite{LS1errata} illustrates this situation.)   
 
\noindent $\left(a_{4}\right)$
$Y$ is conformal on $E$: $Y_{u} \cdot Y_{v} = 0, Y_{u}\cdot Y_{u} = Y_{v}\cdot Y_{v}$
on $E$.
 
\noindent $\left(a_{5}\right)$
$\triangle Y := Y_{uu} + Y_{vv} = H\left(Y\right)Y_{u} \times Y_{v}$  on $E$.
\vspace{1mm} 

\noindent Here by the (open) arcs ${\bf o}_{1}{\bf b}$  and ${\bf o}_{2}{\bf a}$  are meant the 
component of $\partial E\setminus\{{\bf o}_{1},{\bf b}\}$  which does not contain  ${\bf a}$ 
and the component of $\partial E\setminus\{{\bf o}_{2},{\bf a}\}$  which does not contain  ${\bf b}$  
respectively. 
Let  $\sigma_{0} = \partial E \setminus \sigma$.
\vspace{2mm}

There are two cases we will need to consider during the proofs of Theorem \ref{CONCLUSION} and Theorem \ref{DEMONSTRATION}: 
\begin{itemize}
\item[$\left(A\right)$] ${\bf o}_{1}= {\bf o}_{2}.$
\item[$\left(B\right)$] ${\bf o}_{1}\neq {\bf o}_{2}.$
\end{itemize}
These correspond to Cases 5 and 3 respectively in Step 1 of the proof of Theorem~1 of \cite{LS1}. 
\vspace{2mm}

\section{Proof of Theorem \ref{CONCLUSION}}  

\begin{figure}[ht]
\centering
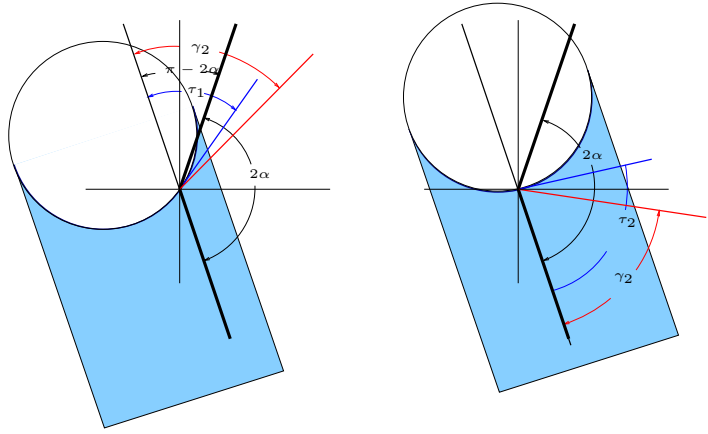
\caption{Left  $\Delta_{\beta_{1}}:$  domain of $h^{-}_{\beta_{1}}$  \ \ \ \  
Right  $\Delta_{\beta_{2}}:$  domain of $h^{+}_{\beta_{2}}$ \label{pizza}}
\end{figure}

Since $\pi-2\alpha<\gamma_{2}<2\alpha,$  we can choose $\tau_{1}\in (\pi-2\alpha, \gamma_{2})$  and 
and $\tau_{2}\in (\gamma_{2},2\alpha).$
Set $\beta_{1}=\frac{\pi}{2}-\tau_{1}$  and  $\beta_{2}=\frac{\pi}{2}-(\pi-\tau_{2})=\tau_{2}-\frac{\pi}{2}.$
With these choices of $\beta_{1}$  and $\beta_{2},$   notice that 
\[
T\left( h^{-}\circ T_{\beta_{1}}\right)(x_{1},0)\cdot (0,-1) = \cos\left(\tau_{1}\right)
> \cos(\gamma_{2}), \ \  {\rm for} \ 0<x_{1}< 2-r_{0}
\]
and 
\[
T\left( h^{+}\circ T_{\beta_{2}}\right)(x_{1},0)\cdot (0,-1) = \cos\left(\tau_{2}\right)
< \cos(\gamma_{2}), \ \  {\rm for} \ 0<x_{1}< 2-r_{0}.
\]
This implies that for $\delta_{1}=\delta_{1}(\beta_{1},\beta_{2})>0$  small enough, 
\begin{equation}
\label{Contact}
T\left( h^{-}_{\beta_{1}}\right)({\bf x})\cdot \vec \nu({\bf x}) > \cos(\gamma({\bf x}))\ \ \ {\rm and} \ \ \ 
T\left( h^{+}_{\beta_{2}}\right)({\bf x})\cdot \vec \nu({\bf x}) < \cos(\gamma({\bf x}))
\end{equation}
for ${\bf x}\in\partial^{-}\Omega$  with $|{\bf x}|<\delta_{1},$  where $\vec \nu({\bf x})$  is the exterior unit normal to 
$\Omega$  at ${\bf x}\in\partial\Omega.$  (See Figure \ref{pizza}.)  
(We may also assume $\nu({\bf x})\cdot (1,1) <0$  for ${\bf x}\in\partial^{+}\Omega$  with $|{\bf x}|<\delta_{1}$   and 
$\nu({\bf x})\cdot (1,-1) <0$  for ${\bf x}\in\partial^{-}\Omega$  with $|{\bf x}|<\delta_{1}$  since $\alpha>\frac{\pi}{4}.$)   

Let $\mu\in (0,\min\{\gamma_{2}-(\pi-2\alpha),2\alpha-\gamma_{2}\})$  and set $\tau_{1}(\mu)=\pi-2\alpha+\mu$  and 
$\tau_{2}(\mu)=2\alpha-\mu,$  so that  $\beta_{1}=\beta_{2}.$  
Let us write $\delta_{1}(\mu)$  for $\delta_{1}(\beta_{1},\beta_{2}),$  $h^{+}_{\mu}$  for $h^{+}_{\beta_{2}},$  
$h^{-}_{\mu}$  for $h^{-}_{\beta_{1}}$  and $\Delta_{\mu}$  for  $\Delta_{\beta_{1}}=\Delta_{\beta_{2}}.$
Since $\beta_{1},\beta_{2}\neq \pm \frac{\pi}{2},$  there exists $R=R(\mu)>0$  such that 
$B({\cal O},R(\mu))\cap \Omega^{*}\subset \Delta_{\mu}$  (where $B({\cal O},R)=\{ {\bf x}\in\Real^{2}:|{\bf x}|<R\}).$

Let us first assume that $\left(A\right)$  holds  and set ${\bf o}={\bf o}_{1}={\bf o}_{2}.$  
\vspace{1mm}

\noindent {\bf Claim:} $f$  is uniformly continuous on $\Omega_{0},$  where $\Omega_{0} = \Omega^{*}\cap \Delta_{\mu}.$
\vspace{1mm}

\noindent {\bf Pf:}  For $r>0,$  set $B_{r}= \{ {\bf u} \in \overline{E} : |{\bf u} - {\bf o}| < r \},$   
$C_{r} = \{ {\bf u} \in \overline{E} : |{\bf u} - {\bf o}| = r \}$  and let $l_{r}$ be the length of the image curve 
$Y(C_{r});$  also let $C^{\prime}_{r} = G(C_{r})$  and $B^{\prime}_{r}= G(B_{r}).$   
From the Courant-Lebesgue Lemma (e.g. Lemma $3.1$ in \cite{Cour:50}), we see that for each $\delta\in (0,1),$ 
there exists a $\rho=\rho(\delta)\in \left(\delta,\sqrt{\delta}\right)$  such that the arclength 
$l_{\rho}$  of $Y(C_{\rho})$  is less than $p(\delta).$   
For $\delta>0,$  let $k(\delta)= \inf_{{\bf u}\in C_{\rho(\delta)}}c({\bf u}) = \inf_{ {\bf x}\in C^{\prime}_{\rho(\delta)} } f({\bf x})$  
and  $m(\delta)= \sup_{{\bf u}\in C_{\rho(\delta)}}c({\bf u}) = \sup_{ {\bf x}\in C^{\prime}_{\rho(\delta)} } f({\bf x});$ 
notice that $m(\delta)-k(\delta)\le l_{\rho} < p(\delta).$  

\begin{figure}[ht]
\centering
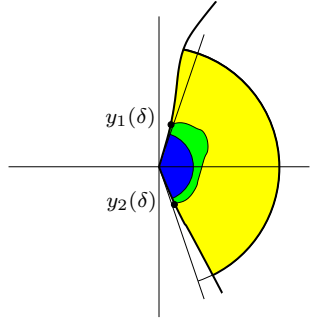
\caption{   $B_{\eta(\delta)}({\cal O})$  (blue region)  and $B^{\prime}_{\rho(\delta)}$  (blue \& green regions)
\label{THREE}}
\end{figure}

For each $\delta\in (0,1)$  with $\sqrt{\delta}<\min\{|{\bf o}-{\bf a}|, |{\bf o}-{\bf b}|\},$  there are two points in 
$C_{\rho(\delta)}\cap\partial E;$   we denote these points as ${\bf e}_{1}(\delta)\in {\bf o}{\bf b}$  and 
${\bf e}_{2}(\delta)\in {\bf o}{\bf a}$  and set ${\bf y}_{1}(\delta)=G({\bf e}_{1}(\delta))$  
and ${\bf y}_{2}(\delta)=G({\bf e}_{2}(\delta)).$
Notice that $C^{\prime}_{\rho(\delta)}$  is a curve in $\overline{\Omega}$  which joins ${\bf y}_{1}\in {\partial}^{+}\Omega^{*}$ 
and ${\bf y}_{2}\in {\partial}^{-}\Omega^{*}$  and 
$\partial\Omega\cap C^{\prime}_{\rho(\delta)}\setminus \{{\bf y}_{1},{\bf y}_{2}\}=\emptyset;$  
therefore there exists $\eta=\eta(\delta)>0$  such that  
$B_{\eta(\delta)}({\cal O})=\{ {\bf x}\in \Omega : |{\bf x}|<\eta(\delta)\} \subset B^{\prime}_{\rho(\delta)}$  
(see Figure \ref{THREE}).

Let $\epsilon>0.$  Choose $\delta>0$  such that $\sqrt{\delta}<\min\{|{\bf o}-{\bf a}|, |{\bf o}-{\bf b}|\},$  
$p(\delta)<\delta_{1}(\mu),$   $p(\delta)<R(\mu),$   and $p(\delta)+q(p(\delta))<\frac{1}{2}\epsilon.$  
Pick a point ${\bf w}\in C^{\prime}_{\rho(\delta)}$  and define $b^{\pm}_{j}:\Delta_{\mu}\to\Real$  by 
\[
b^{\pm}({\bf x})=f({\bf w})  \pm p(\delta)  \pm h_{\mu}^{\mp}({\bf x}), \ \ \ \ {\bf x}\in \Delta_{\mu}.
\]
From (\ref{MeanBounds}), (\ref{Contact}) and the General Comparison Principle (e.g. \cite{FinnBook}, Theorem 5.1),  
we have 
\[
b^{-}({\bf x})<f({\bf x})<b^{+}({\bf x}) \ \ \ \ {\rm for \ all} \ \  {\bf x}\in B^{\prime}_{\rho(\delta)}  \cap \Delta_{\mu}.
\]
Thus if ${\bf x}_{1},{\bf x}_{2}\in \Omega_{0}$  satisfy $|{\bf x}_{1}|<\eta(\delta),$  $|{\bf x}_{2}|<\eta(\delta)$  and 
$|{\bf x}_{1}-{\bf x}_{2}|<\eta(\delta),$   then 
\begin{equation}
\label{Peach}
|f({\bf x}_{1})-f({\bf x}_{2})| <2p(\delta)+2q\left(p(\delta)\right)<\epsilon.
\end{equation}

Since $f$  is uniformly continuous on $\{ {\bf x}\in \Omega^{*} \ : \ |{\bf x}|\ge \frac{1}{2}\eta(\delta)\},$  
there exists a $\lambda>0$  such that if ${\bf x}_{1},{\bf x}_{2}\in \Omega^{*}$  satisfy 
$|{\bf x}_{1}|\ge \frac{1}{2}\eta(\delta),$  ${\bf x}_{2}|\ge \frac{1}{2}\eta(\delta)$  and 
$|{\bf x}_{1}-{\bf x}_{2}| < \lambda,$  then  $|f({\bf x}_{1})-f({\bf x}_{2})|<\epsilon.$
Now set $d=d(\epsilon)=\min\{\lambda, \frac{1}{2}\eta(\delta)\}.$  
If ${\bf x}_{1},{\bf x}_{2}\in \Omega_{0},$  $|{\bf x}_{1}-{\bf x}_{2}|< d(\epsilon)\le \frac{1}{2}\eta(\delta)$  
and $|{\bf x}_{1}|<\frac{1}{2}\eta(\delta),$  then $|{\bf x}_{1}|<\eta(\delta)$  and $|{\bf x}_{2}|<\eta(\delta);$  
hence $|f({\bf x}_{1})-f({\bf x}_{2})|<\epsilon$  by (\ref{Peach}).    
Next, if ${\bf x}_{1},{\bf x}_{2}\in \Omega_{0},$  $|{\bf x}_{1}-{\bf x}_{2}|< d(\epsilon)\le\lambda,$ 
$|{\bf x}_{1}|\ge \frac{1}{2}\eta(\delta)$  and $|{\bf x}_{2}|\ge \frac{1}{2}\eta(\delta),$
then $|f({\bf x}_{1})-f({\bf x}_{2})|<\epsilon.$  
Therefore, for all ${\bf x}_{1},{\bf x}_{2}\in \Omega_{0}$  with $|{\bf x}_{1}-{\bf x}_{2}|< d(\epsilon),$  
we have $|f({\bf x}_{1})-f({\bf x}_{2})|<\epsilon.$   The claim is proven.
\vspace{2mm}

Notice that if $\theta(\mu)= \alpha-\mu$  ($=\tau_{2}(\mu)-\alpha =\pi-\alpha-\tau_{1}(\mu)$), then 
\[
\{\left(r\cos(\theta(\mu)),r\sin(\theta(\mu))\right) : r\ge 0\}
\]
is the tangent ray to $\partial\Omega_{0}$  at ${\cal O}$  and it follows from the Claim that 
$f\in C^{0}(\overline{\Omega_{0}});$  hence the radial limits $Rf(\theta)$  of $f$  at ${\cal O}$  exist for 
$\theta\in [-\alpha,\theta(\mu)]$  and the radial limits are identical 
(i.e. $Rf(\theta)=f({\cal O})$  for all $\theta\in [-\alpha,\theta(\mu)],$  where $f({\cal O})$  is the value at ${\cal O}$  
of the restriction of $f$  to $\overline{\Omega_{0}}$).  
Since $\lim_{\mu\downarrow 0} \theta(\mu) = \alpha,$  Theorem \ref{CONCLUSION} is proven in this case.
\begin{figure}[ht]
\centering
\includegraphics{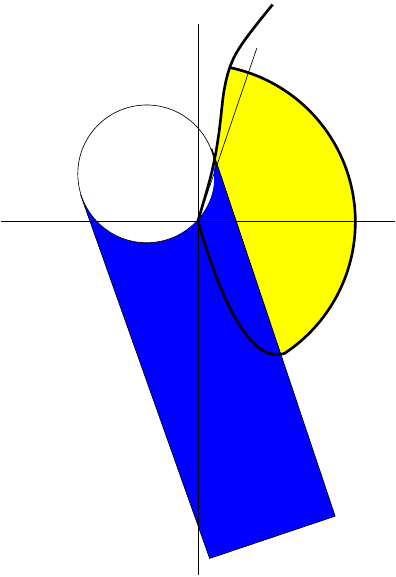}
\caption{The domain (in blue) of the toroidal functions $h^{\pm}_{\mu},$  $\alpha>\frac{\pi}{4}.$    \label{ONEa1}}
\end{figure}
\vspace{2mm}

Let us next assume that $\left(B\right)$ holds.  This part of the proof is essentially the same as the proof of case $\left(B\right)$  
in Theorem 1 of \cite{NoraKirk1}.  As in  \cite{NoraKirk1} and taking the hypothesis $\alpha\le \frac{\pi}{2}$  into account, we see that  
\begin{itemize}
\item[(i)] $c\in C^{0}\left(\overline{E}\setminus \{{\bf o}_{1}, {\bf o}_{2} \} \right),$  
\item[(ii)] there exist $\alpha_{1},\alpha_{2}\in [-\alpha,\alpha]$  with $\alpha_{1}<\alpha_{2}$  such that 
$Rf(\theta)$  exists when $\theta\in \left(\alpha_{1},\alpha_{2}\right),$    and 
\item[(iii)] $Rf$  is strictly increasing or  strictly  decreasing  on $(\alpha_{1}, \alpha_{2}).$ 
\end{itemize}
Taking hypothesis (\ref{Limit-}) into account and using cylinders as in Case 3 of Step 1 in the proof of Theorem 1 of 
\cite{LS1} (see Figure 2b in \cite{LS1errata}) or using $h^{\pm}_{\mu}$  (see Figure \ref{ONEa1}), we see that in addition to 
(i)-(iii), we have 
\begin{itemize}
\item[(iv)] $c\in C^{0}\left(\overline{E}\setminus \{{\bf o}_{1} \} \right)$  and 
\item[(v)] $Rf(\theta)$  exists when $\theta\in \left[-\alpha,\alpha_{2}\right).$
\end{itemize}
If $\alpha_{2}=\alpha,$  then Theorem \ref{CONCLUSION} is proven.  
Otherwise, suppose $\alpha_{2}<\alpha$  and fix $\delta_{0}\in (0,\delta^{*})$  and $\Omega_{0} = \Omega^{*}\cap \Delta_{\mu}$
 as before.  
\vspace{2mm}

\noindent {\bf Claim:} Suppose $\alpha_{2}<\alpha.$  Then $f$  is uniformly continuous on $\Omega^{+}_{0},$  where 
\[
\Omega^{+}_{0} \myeq \{(r\cos(\theta),r\sin(\theta))\in \Omega_{0} : 0<r<\delta^{*}, \alpha_{2}<\theta<\pi\}.
\]
\vspace{2mm}

\noindent Notice that the restriction of $Y$  to $G^{-1}\left(\overline{\Omega^{+}_{0}}\right)$  maps only one point, ${\bf o}_{1},$  
to ${\cal O}\times\Real$  and so the proof of this claim is the same as the proof of the previous Claim.
Thus $f\in C^{0}\left(\overline{\Omega^{+}_{0}}\right);$  since $\lim_{\mu\downarrow 0} \theta(\mu) = \alpha,$ 
we see that 
\[
Rf(\theta)=\lim_{\tau\uparrow \alpha_{2}} Rf(\tau) \ \ \ \ {\rm for \ all} \ \ \theta\in [\alpha_{2},\alpha).
\]  
Thus Theorem \ref{CONCLUSION} is proven.

\section{Proof of Theorem \ref{DEMONSTRATION}}  

\begin{figure}[ht]
\centering
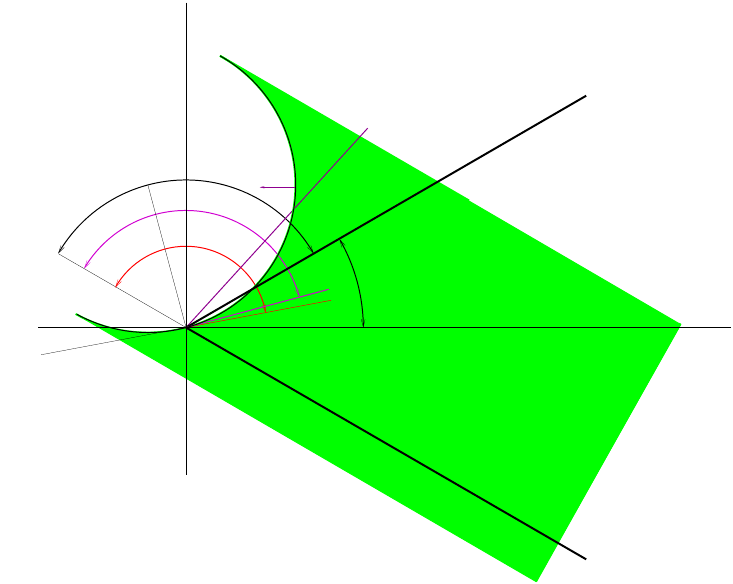
\caption{$\alpha=\frac{\pi}{6},$  $\lambda_{1}=0,$  $\lambda_{2}=\frac{\pi}{2},$  $\gamma_{2}=\frac{7\pi}{9},$  and 
$\tau_{1}=\frac{27\pi}{36}.$  
The domain of $h^{-}_{\beta_{1}}$  is the green region.  \label{ZEROX2}}
\end{figure}

\begin{figure}[ht]
\centering
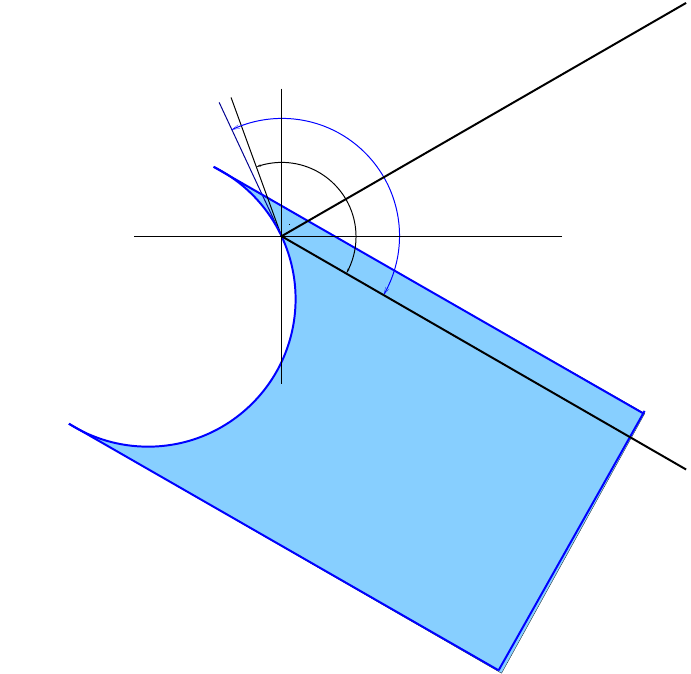
\caption{$\alpha=\frac{\pi}{6},$  $\lambda_{1}=0,$  $\lambda_{2}=\frac{\pi}{2},$  $\gamma_{2}=\frac{7\pi}{9},$  and 
$\tau_{2}=\frac{29\pi}{36}.$  
The domain of $h^{+}_{\beta_{2}}$  is the blue region.  \label{ZEROX3}}
\end{figure}

Suppose (\ref{Condition-1}) does not hold.  
Since $\pi-2\alpha-\lambda_{1}<\gamma_{2}<\pi+2\alpha-\lambda_{2},$  we can choose $\tau_{1},\tau_{2}\in (0,\pi)$  
such that $\tau_{1}\in (\pi-2\alpha-\lambda_{1}, \gamma_{2})$   and $\tau_{2}\in (\gamma_{2},\pi+2\alpha-\lambda_{2}).$
Set $\beta_{1}=\frac{\pi}{2}-\tau_{1}$  and  $\beta_{2}=\tau_{2}-\frac{\pi}{2}.$  
(See Figures \ref{ZEROX2} and \ref{ZEROX3}.)  
With these choices of $\beta_{1}$  and $\beta_{2},$   notice that 
\[
T\left( h^{-}\circ T_{\beta_{1}}\right)(x_{1},0)\cdot (0,-1) = \cos\left(\tau_{1}\right)
> \cos(\gamma_{2}), \ \  {\rm for} \ 0<x_{1}< 2-r_{0}
\]
and 
\[
T\left( h^{+}\circ T_{\beta_{2}}\right)(x_{1},0)\cdot (0,-1) = \cos\left(\tau_{2}\right)
< \cos(\gamma_{2}), \ \  {\rm for} \ 0<x_{1}< 2-r_{0}.
\]
This implies that for $\delta_{1}=\delta_{1}(\beta_{1},\beta_{2})>0$  small enough, 
\begin{equation}
\label{ContactA}
T\left( h^{-}_{\beta_{1}}\right)({\bf x})\cdot \vec \nu({\bf x}) > \cos(\gamma({\bf x}))\ \ \ {\rm and} \ \ \ 
T\left( h^{+}_{\beta_{2}}\right)({\bf x})\cdot \vec \nu({\bf x}) < \cos(\gamma({\bf x}))
\end{equation}
for ${\bf x}\in\partial^{-}\Omega$  with $|{\bf x}|<\delta_{1},$  where $\vec \nu({\bf x})$  is the exterior unit normal to 
$\Omega$  at ${\bf x}\in\partial\Omega.$  (See Figures \ref{pizza}, \ref{ZEROX2} and \ref{ZEROX3}.)

Notice that the tangent plane at $(0,0,0)$  to the surface 
$\{({\bf x}, h^{-}_{\beta_{1}}({\bf x})) : {\bf x}\in \Delta_{\beta_{1}}\}$  
is a vertical plane with (downward oriented) unit normal $\vec n=(-\sin(\tau_{1}+\alpha),-\cos(\tau_{1}+\alpha),0)$  and   
$\lim_{\partial^{+}\Omega\ni {\bf x}\to {\cal O}} \vec\nu({\bf x})=(-\sin(\alpha),\cos(\alpha),0).$  
Suppose $\tau_{1}+2\alpha\le\pi.$  Then 
\[
\lim_{\partial^{+}\Omega\ni {\bf x}\to {\cal O}} \vec n \cdot \vec\nu({\bf x})=  -\cos(\tau_{1}+2\alpha) > 
-\cos(\pi-\lambda_{1}) = \cos(\lambda_{1})
\]
since $\tau_{1}+2\alpha>\pi-\lambda_{1};$  since 
$\liminf_{\partial^{+}\Omega\ni {\bf x}\to {\cal O}} \gamma({\bf x})\ge \lambda_{1},$
this implies that for some $\delta_{2}>0$  small enough, 
\begin{equation}
\label{Contact2}
T\left( h^{-}_{\beta_{1}}\right)({\bf x})\cdot \vec \nu({\bf x}) > \cos(\gamma({\bf x})) \ \ \ \ 
{\rm for} \ \ {\bf x}\in\partial^{+}\Omega \ \ {\rm with} \ \ |{\bf x}|<\delta_{2}.
\end{equation}
If $\tau_{1}+2\alpha>\pi,$  then $\lambda_{1}$  doesn't matter and we argue as in the proof of Theorem \ref{CONCLUSION}; 
see Figure \ref{ZEROX2} for an illustration of this case.

Now the tangent plane at $(0,0,0)$  to the surface 
$\{({\bf x}, h^{+}_{\beta_{2}}({\bf x})) : {\bf x}\in \Delta_{\beta_{2}}\}$  
is a vertical plane with (downward oriented) unit normal $\vec m=(\sin(\tau_{2}-\alpha),-\cos(\tau_{2}-\alpha),0)$  and   
$\lim_{\partial^{+}\Omega\ni {\bf x}\to {\cal O}} \vec\nu({\bf x})=(-\sin(\alpha),\cos(\alpha),0).$  
Suppose $\tau_{2}\ge 2\alpha.$  Then 
\[
\lim_{\partial^{+}\Omega\ni {\bf x}\to {\cal O}} \vec m \cdot \vec\nu({\bf x})=  -\cos(\tau_{2}-2\alpha) <  
-\cos(\pi-\lambda_{2}) = \cos(\lambda_{2})
\]
since $\tau_{2}-2\alpha<\pi-\lambda_{2};$  since 
$\limsup_{\partial^{+}\Omega\ni {\bf x}\to {\cal O}} \gamma({\bf x})\le \lambda_{2},$
this implies that for some $\delta_{3}>0$  small enough, 
\begin{equation}
\label{Contact4}
T\left( h^{+}_{\beta_{1}}\right)({\bf x})\cdot \vec \nu({\bf x}) < \cos(\gamma({\bf x})) \ \ \ \ 
{\rm for} \ \ {\bf x}\in\partial^{+}\Omega \ \ {\rm with} \ \ |{\bf x}|<\delta_{3}.
\end{equation}
If $\tau_{2}<2\alpha,$  then $\lambda_{2}$  doesn't matter and we argue as in the proof of Theorem \ref{CONCLUSION}. 

Now set $\delta_{4}=\min \{ \delta_{1}, \delta_{2}, \delta_{3} \}.$
The proof of Theorem \ref{DEMONSTRATION} now follows essentially as in the proof of Theorem \ref{CONCLUSION}.

\end{document}

%% file: figure_2_12a.pdf_t
\begin{picture}(0,0)%
\includegraphics{figure_2_12a.pdf}%
\end{picture}%
\setlength{\unitlength}{829sp}%
\begingroup\makeatletter\ifx\SetFigFont\undefined%
\gdef\SetFigFont#1#2#3#4#5{%
  \reset@font\fontsize{#1}{#2pt}%
  \fontfamily{#3}\fontseries{#4}\fontshape{#5}%
  \selectfont}%
\fi\endgroup%
\begin{picture}(22524,9024)(-4511,-3673)
\put(-1197,1233){\makebox(0,0)[lb]{\smash{{\SetFigFont{6}{7.2}{\rmdefault}{\mddefault}{\updefault}{$\tau_{1}$}%
}}}}
\put(14521,1319){\makebox(0,0)[lb]{\smash{{\SetFigFont{6}{7.2}{\rmdefault}{\mddefault}{\updefault}{$\tau_{2}$}%
}}}}
\put(14386,419){\makebox(0,0)[lb]{\smash{{\SetFigFont{6}{7.2}{\rmdefault}{\mddefault}{\updefault}{$-\beta$}%
}}}}
\put(1021,1068){\makebox(0,0)[lb]{\smash{{\SetFigFont{6}{7.2}{\rmdefault}{\mddefault}{\updefault}{$\beta$}%
}}}}
\end{picture}%

%% file: figure_2_9E.pdf_t
\begin{picture}(0,0)%
\includegraphics{figure_2_9E.pdf}%
\end{picture}%
\setlength{\unitlength}{1036sp}%
\begingroup\makeatletter\ifx\SetFigFont\undefined%
\gdef\SetFigFont#1#2#3#4#5{%
  \reset@font\fontsize{#1}{#2pt}%
  \fontfamily{#3}\fontseries{#4}\fontshape{#5}%
  \selectfont}%
\fi\endgroup%
\begin{picture}(16731,10189)(1303,-8038)
\put(5671,1019){\makebox(0,0)[lb]{\smash{{\SetFigFont{5}{6.0}{\rmdefault}{\mddefault}{\updefault}{$\gamma_{2}$}%
}}}}
\put(5041,479){\makebox(0,0)[lb]{\smash{{\SetFigFont{5}{6.0}{\rmdefault}{\mddefault}{\updefault}{$\pi-2\alpha$}%
}}}}
\put(15796,-4471){\makebox(0,0)[lb]{\smash{{\SetFigFont{5}{6.0}{\rmdefault}{\mddefault}{\updefault}{$\gamma_{2}$}%
}}}}
\put(7072,-2035){\makebox(0,0)[lb]{\smash{{\SetFigFont{5}{6.0}{\rmdefault}{\mddefault}{\updefault}{$2\alpha$}%
}}}}
\put(15886,-3121){\makebox(0,0)[lb]{\smash{{\SetFigFont{5}{6.0}{\rmdefault}{\mddefault}{\updefault}{$\tau_{2}$}%
}}}}
\put(5581,-16){\makebox(0,0)[lb]{\smash{{\SetFigFont{5}{6.0}{\rmdefault}{\mddefault}{\updefault}{$\tau_{1}$}%
}}}}
\put(15031,-1591){\makebox(0,0)[lb]{\smash{{\SetFigFont{5}{6.0}{\rmdefault}{\mddefault}{\updefault}{$2\alpha$}%
}}}}
\end{picture}%

%% file: figure_2_2b.pdf_t
\begin{picture}(0,0)%
\includegraphics{figure_2_2b.pdf}%
\end{picture}%
\setlength{\unitlength}{1657sp}%
\begingroup\makeatletter\ifx\SetFigFont\undefined%
\gdef\SetFigFont#1#2#3#4#5{%
  \reset@font\fontsize{#1}{#2pt}%
  \fontfamily{#3}\fontseries{#4}\fontshape{#5}%
  \selectfont}%
\fi\endgroup%
\begin{picture}(4524,4770)(3139,-4573)
\put(4616,-2916){\makebox(0,0)[lb]{\smash{{\SetFigFont{8}{9.6}{\rmdefault}{\mddefault}{\updefault}{$y_{2}(\delta)$}%
}}}}
\put(4601,-1651){\makebox(0,0)[lb]{\smash{{\SetFigFont{8}{9.6}{\rmdefault}{\mddefault}{\updefault}{$y_{1}(\delta)$}%
}}}}
\end{picture}%

%% file: figure_2_11f.pdf_t
\begin{picture}(0,0)%
\includegraphics{figure_2_11f.pdf}%
\end{picture}%
\setlength{\unitlength}{622sp}%
\begingroup\makeatletter\ifx\SetFigFont\undefined%
\gdef\SetFigFont#1#2#3#4#5{%
  \reset@font\fontsize{#1}{#2pt}%
  \fontfamily{#3}\fontseries{#4}\fontshape{#5}%
  \selectfont}%
\fi\endgroup%
\begin{picture}(22290,17706)(1078,-12355)
\put(7046,-63){\makebox(0,0)[lb]{\smash{{\SetFigFont{6}{7.2}{\rmdefault}{\mddefault}{\updefault}{$\pi-2\alpha$}%
}}}}
\put(7201,-1051){\makebox(0,0)[lb]{\smash{{\SetFigFont{6}{7.2}{\rmdefault}{\mddefault}{\updefault}{$\tau_{1}$}%
}}}}
\put(7381,-2266){\makebox(0,0)[lb]{\smash{{\SetFigFont{6}{7.2}{\rmdefault}{\mddefault}{\updefault}{$\gamma_{2}$}%
}}}}
\put(12086,-3198){\makebox(0,0)[lb]{\smash{{\SetFigFont{6}{7.2}{\rmdefault}{\mddefault}{\updefault}{$\alpha$}%
}}}}
\end{picture}%

%% file: figure_2_11g.pdf_t
\begin{picture}(0,0)%
\includegraphics{figure_2_11g.pdf}%
\end{picture}%
\setlength{\unitlength}{622sp}%
\begingroup\makeatletter\ifx\SetFigFont\undefined%
\gdef\SetFigFont#1#2#3#4#5{%
  \reset@font\fontsize{#1}{#2pt}%
  \fontfamily{#3}\fontseries{#4}\fontshape{#5}%
  \selectfont}%
\fi\endgroup%
\begin{picture}(20970,20496)(-1835,-17893)
\put(7246,-2491){\makebox(0,0)[lb]{\smash{{\SetFigFont{6}{7.2}{\rmdefault}{\mddefault}{\updefault}{$\gamma_{2}$}%
}}}}
\put(7336,-1051){\makebox(0,0)[lb]{\smash{{\SetFigFont{6}{7.2}{\rmdefault}{\mddefault}{\updefault}{$\tau_{2}$}%
}}}}
\end{picture}%